\documentclass[12pt]{amsart}
\usepackage{amssymb}
\usepackage{graphicx}

\theoremstyle{plain}
\newtheorem{theorem}{Theorem}[section]
\newtheorem{lemma}[theorem]{Lemma}
\newtheorem{proposition}[theorem]{Proposition}
\newtheorem{corollary}[theorem]{Corollary}
\theoremstyle{definition}

\newtheorem{remark}[theorem]{Remark}

\numberwithin{equation}{section}

\renewcommand{\phi}{\varphi}
\renewcommand{\theta}{\vartheta}
\newcommand{\R}{\mathbb{R}}
\newcommand{\Z}{\mathbb{Z}}
\newcommand{\T}{\mathbb{T}}
\newcommand{\C}{\mathbb{C}}
\newcommand{\zero}{{\mathbf{0}}}
\newcommand{\eps}{\varepsilon}
\newcommand{\aand}{\textrm{\ \ and\ \ }}
\newcommand{\vecd}{{\mathbf{d}}}
\newcommand{\veci}{{\mathbf{i}}}
\newcommand{\vecj}{{\mathbf{j}}}
\newcommand{\veck}{{\mathbf{k}}}
\newcommand{\vecl}{{\mathbf{l}}}
\newcommand{\vecp}{{\mathbf{p}}}
\newcommand{\vecu}{{\mathbf{u}}}
\newcommand{\vecv}{{\mathbf{v}}}
\newcommand{\vecw}{{\mathbf{w}}}
\newcommand{\vecx}{{\mathbf{x}}}
\newcommand{\vecy}{{\mathbf{y}}}
\newcommand{\vecz}{{\mathbf{z}}}
\newcommand{\vecalpha}{{\mathbf{\alpha}}}

\begin{document}

\title{Rotation sets of billiards with one obstacle}
\author{Alexander Blokh}
\address{Department of Mathematics, University of Alabama in
  Birmingham, University Station, Birmingham, AL  35294-2060}
\email{ablokh@math.uab.edu}
\author{Micha\l\ Misiurewicz}
\address{Department of Mathematical Sciences, IUPUI, 402 N. Blackford
  Street, Indianapolis, IN 46202-3216}
\email{mmisiure@math.iupui.edu}
\author{N\'andor Sim\'anyi}
\address{Department of Mathematics, University of Alabama in
  Birmingham, University Station, Birmingham, AL  35294-2060}
\email{simanyi@math.uab.edu}
\thanks{The first author was partially supported by NSF grant DMS
   0456748, the second author by NSF grant DMS 0456526, and the third
   author by NSF grant DMS 0457168}
\date{August 16, 2005}
\maketitle
\begin{abstract}
   We investigate the rotation sets of billiards on the
   $m$-dimensional torus with one small convex obstacle and in the
   square with one small convex obstacle. In the first case the
   displacement function, whose averages we consider, measures the change of the
   position of a point in the universal
   covering of the torus (that is, in the Euclidean space), in the
   second case it measures the rotation around the obstacle. A substantial
   part of the rotation set has usual strong properties of rotation
   sets.
\end{abstract}

\section{Introduction}\label{sec_intro}

Traditionally, billiards have been investigated from the point of
view of Ergodic Theory. That is, the properties that have been
studied, were the statistical properties with respect to the natural
invariant measure equivalent to the Lebesgue measure. However, it is
equally important to investigate the limit behavior of \emph{all}
trajectories and not only of \emph{almost all} of them. In
particular, periodic trajectories (which are of zero measure) are of
great interest. A widely used method in this context is to observe
that billiards in convex domains are twist maps (see, e.g.\
\cite{KH}, Section~9.2), so well developed rotation theory for twist
maps (see, e.g.\ \cite{KH}, Section~9.3) applies to them.

Rotation Theory has been recently developed further and its scope
has been significantly widened (see, e.g., Chapter~6 of \cite{ALM}
for a brief overview). This opens possibilities of its application
to new classes of billiards. In the general Rotation Theory one
considers a dynamical system together with an \emph{observable},
that is a function on the phase space, with values in a vector
space. Then one takes limits of ergodic averages of the observable
along longer and longer pieces of trajectories. The \emph{rotation
set} obtained in such a way contains all averages of the observable
along periodic orbits, and, by the Birkhoff Ergodic Theorem,
integrals of the observable with respect to all ergodic invariant
probability measures. With the natural choice of an observable,
information about the appropriate rotation set allows one to
describe the behavior of the trajectories of the system (see examples in
Chapter~6 of \cite{ALM}). Exact definitions are given later in the
paper.

We have to stress again that we consider \emph{all} trajectories.
Indeed, restricting attention to one ergodic measure would result in
seeing only one rotation vector. However, rotation vectors of other
points, non-typical for the measure, will be missing. Thus, the
approach to the billiards should be from the point of view of
Topological Dynamics, instead of Ergodic Theory, even though we do
consider various invariant measures for which rotation vector can be
computed. Observe that the rotation set for a suitably chosen
observable is a useful characteristic of the dynamical system.

In the simplest case, the observable is the increment in one step
(for discrete systems) or the derivative (for systems with
continuous time) of another function, called \emph{displacement}.
When the displacement is chosen in a natural way, the results on the
rotation set are especially interesting. In this paper we consider
two similar classes of billiards, and the observables which we use
for them are exactly of that type. One system consists of billiards
on an $m$-dimensional torus with one small convex obstacle which we
lift to the universal covering of the torus (that is, to the
Euclidean space) and consider the natural displacement there.
Those  models constitute the rigorous mathematical formulation of the
so called Lorentz gas dynamics with periodic configuration of
obstacles. They are especially important for physicists doing
research in the foundations of nonequilibrium dynamics, since the
Lorentz gas serves as a good paradigm for nonequilibrium stationary
states, see the nice survey \cite{Dettmann}.

The other system consists of billiards in a square with one small convex
obstacle close to the center of the square; here we measure average
rotation around a chosen obstacle using the argument as the
displacement.

We treat both billiards as flows. This is caused by the fact that in
the lifting (or unfolding for a billiard in a square) we may have
infinite horizon, especially if the obstacle is small. In other
words, there are infinite trajectories without reflections, so when
considering billiards as maps, we would have to divide by zero.
Although this is not so bad by itself (infinity exists), we lose
compactness and cannot apply nice general machinery of the rotation
theory (see, e.g., \cite{Z}).

Note that in the case of a billiard in a square we have to deal with
trajectories that reflect from the vertices of the square. We can
think about such reflection as two infinitesimally close
reflections from two adjacent sides. Then it is clear that our
trajectory simply comes back along the same line on which it arrived to the
vertex, and that this does not destroy the continuity of the flow.

The ideas, methods, and results in both cases, the torus and the
square, are very similar. However, there are some important
differences, and, in spite of its two-dimensionali\-ty, in general the
square case is more complicated. Therefore we decided to treat the
torus case first (in Sections~\ref{sec_prelim}, \ref{sec_rotset}
and~\ref{sec_arlarge}) and then, when we describe the square case
(in Section~\ref{sec_square}), we describe the differences from the
torus case, without repeating the whole proofs. We believe that this
type of exposition is simpler for a reader than the one that treats
both cases simultaneously or the one that produces complicated
abstract theorems that are then applied in both cases. In
Section~\ref{sec_conn} we get additionally some general results,
applicable also to other situations.

Let us describe shortly the main results of the paper. The exact
definitions will be given later. Let us only note that the
\emph{admissible rotation set} is a subset of the full rotation set,
about which we can prove much stronger results than about the full
rotation set. Also, a \emph{small} obstacle does not mean
``arbitrarily small'' one. We derive various estimates of the size
of the admissible rotation set. In the torus case the estimates that
are independent of the dimension are non-trivial because of the
behavior of the geometry of $\R^m$ as $m\to\infty$. In both cases we
show that the admissible rotation set approximates better and better
the full rotation set when the size of the obstacle diminishes.

We prove that in both cases, the torus and the square, if the
obstacle is small, then the admissible rotation set is convex,
rotation vectors of periodic orbits are dense in it, and if $u$ is a
vector from its interior, then there exists a trajectory with
rotation vector $u$ (and even an ergodic invariant measure, for
which the integral of the velocity is equal to $u$, so that $u$ is
the rotation vector of almost every trajectory). The full rotation
set is connected, and in the case of the square, is equal to the
interval $[-\sqrt{2}/4,\sqrt{2}/4]$.

We conjecture that the full rotation set shares the strong
properties of the admissible rotation set.

\section{Preliminary results - torus}\label{sec_prelim}

Let us consider a billiard on the $m$-dimensional torus
$\T^m=\R^m/\Z^m$ ($m\ge 2$) with one strictly convex (that is, it is
convex and its boundary does not contain any straight line segment)
obstacle $O$ with a smooth boundary.
We do not specify explicitly how large the obstacle is, but let us
think about it as a rather small one. When we lift the whole picture
to $\R^m$ then we get a family of obstacles $O_\veck$, where $\veck\in\Z^m$
and $O_\veck$ is $O_\zero$ translated by the vector $\veck$.

When we speak about a
\emph{trajectory}, we mean a positive (one-sided) billiard
trajectory, unless we explicitly say that it is a full (two-sided)
one. However, we may mean a trajectory in the phase space, in the
configuration space (on the torus), or in the lifting or unfolding
(the Euclidean space). It will be usually clear from the context,
which case we consider.

We will say that the obstacle $O_\veck$ is \emph{between} $O_\veci$ and $O_\vecj$
if it intersects the convex hull of $O_\veci\cup O_\vecj$ and $\veck\ne\veci,\vecj$. For
a trajectory $P$, beginning on a boundary of $O$,
its \emph{type} is a sequence $(\veck_n)_{n=0}^\infty$ of elements of
$\Z^m$ if the continuous lifting of $P$ to $\R^m$ that starts at the
boundary of $O_{\veck_0}$ reflects consecutively from $O_{\veck_n}$,
$n=1,2,\dots$. In order to make the type unique for a given $P$, we
will additionally assume that $\veck_0=\zero$. Note that except for the case
when $P$ at its initial point is tangent to $O$, there are
infinitely many reflections, so the type of $P$ is well defined. This
follows from the following lemma. In it, we do not count tangency as a
reflection.

\begin{lemma}\label{infref}
If a trajectory has one reflection then it has infinitely
many reflections.
\end{lemma}

\begin{proof}
Suppose that a trajectory has a reflection, but there are only
finitely many of them. Then we can start the trajectory from the last
reflection. Its $\omega$-limit set is an affine subtorus of $\T^m$ and
the whole (positive) trajectory is contained in this subtorus (and it
is dense there). Since
we started with a reflection, this subtorus intersects the interior of
the obstacle. Since the trajectory is dense in this subtorus, we get a
contradiction.
\end{proof}

Of course, there may be trajectories without any reflections. In
particular, if the obstacle is contained in a ball of radius less than
$1/2$, there are such trajectories in the direction of the basic unit
vectors.

Sometimes we will speak about the \emph{type} of a piece of a
trajectory; then it is a finite sequence. We will also use the term
\emph{itinerary}.

We will call a sequence $(\veck_n)_{n=0}^\infty$ of elements of $\Z^m$
\emph{admissible} if
\begin{enumerate}
\item\label{adm1} $\veck_0=\zero$,
\item\label{adm2} for every $n$ we have $\veck_{n+1}\ne \veck_n$,
\item\label{adm3} for every $n$ there is no obstacle between $O_{\veck_n}$ and
  $O_{\veck_{n+1}}$,
\item\label{adm4} for every $n$ the obstacle $O_{\veck_{n+1}}$ is not between
  $O_{\veck_n}$ and $O_{\veck_{n+2}}$.
\end{enumerate}

\begin{theorem}\label{thmadm}
For any admissible sequence $(\veck_n)_{n=0}^\infty$ there is a trajectory
with type $(\veck_n)_{n=0}^\infty$. If additionally there is $\vecp\in\Z^m$
and a positive integer $q$ such that $\veck_{n+q}=\veck_n+\vecp$ for every $n$
then this trajectory can be chosen periodic of discrete period $q$
(that is, after $q$ reflections we come back to the starting point in
the phase space). Similarly, for any admissible sequence
$(\veck_n)_{n=-\infty}^\infty$ there is a trajectory with type
$(\veck_n)_{n=-\infty}^\infty$.
\end{theorem}

\begin{proof}
Fix $n$. For every sequence $A=(\vecx_i)_{i=0}^n$ such that $\vecx_i$ belongs
to the boundary of $O_{\veck_i}$ for $i=0,1,\dots,n$, let $\Gamma(A)$ be
the curve obtained by joining consecutive points $\vecx_i$ by straight
segments (such a curve may intersect interiors of some obstacles).
Since the Cartesian product of the boundaries of $O_{\veck_i}$
is compact and the length of $\Gamma(A)$ depends continuously on $A$,
there is an $A$ for which this length is minimal.

We claim that in such a case $\Gamma(A)$ is a piece of a
trajectory. By (\ref{adm3}), the segment $I_i$ joining $\vecx_i$ with
$\vecx_{i+1}$ cannot intersect any obstacle except $O_{\veck_i}$ and
$O_{\veck_{i+1}}$. If it intersects $O_{\veck_i}$ at more than one point, it
intersects its boundary at $\vecx_i$ and at another point $\vecy$. Then
replacing $\vecx_i$ by $y$ will make $\Gamma(A)$ shorter, a contradiction.
This argument does not work only if $\vecx_{i-1}$, $\vecx_i$ and $\vecx_{i+1}$ are
collinear and $\vecx_i$ lies between $\vecx_{i-1}$ and $\vecx_{i+1}$. However,
such situation is excluded by (\ref{adm4}). This proves that $I_i$
does not intersect $O_{\veck_i}$ at more than one point. Similarly, it
does not intersect $O_{\veck_{i+1}}$ at more than one point. Now the known
property of curves with minimal lengths guarantees that at every
$\vecx_i$, $i=1,2,\dots,n-1$, the incidence and reflection angles are
equal. This proves our claim. For a two-sided sequence
$(\veck_n)_{n=-\infty}^\infty$ the argument is very similar.

Now we make this construction for every $n$ and get a sequence
$(A_n)_{n=1}^\infty$ of pieces of trajectories. We note their
initial points in the phase space (points and directions) and choose a
convergent subsequence of those. Then the trajectory of this limit
point in the phase space will have the prescribed type.

If there is $\vecp\in\Z^m$ and a positive integer $q$ such that
$\veck_{n+q}=\veck_n+\vecp$ for every $n$, then we consider only the sequence
$A=(\vecx_i)_{i=0}^{q-1}$ and repeat the first part of the above proof
adding the segment joining $\vecx_{q-1}$ with $\vecx_0+\vecp$ to $\Gamma(A)$.
\end{proof}

Note that by Corollary~1.2 of \cite{Ch82}, if the obstacle is strictly
convex then a periodic orbit from the above theorem is unique.

The next lemma essentially expresses the fact that any billiard flow
with convex obstacles lacks focal points. It follows from the corollary
after Lemma~2 of \cite{St89}. The types of trajectory pieces about which we
speak in this lemma are not necessarily admissible.

\begin{lemma}\label{unique}
For a given finite sequence $B=(\veck_n)_{n=0}^s$ of elements of $\Z^m$
and points $\vecx_0,\vecx_s$ on the boundaries of $O_{\veck_0}$ and $O_{\veck_s}$
respectively, there is at most one trajectory piece of type $B$
starting at $\vecx_0$ and ending at $\vecx_s$. The same remains true if we
allow the first segment of the trajectory piece to cross $O_{\veck_0}$ and
the last one to cross $O_{\veck_s}$ (as in Figure~\ref{cross}).
\end{lemma}

\begin{figure}\refstepcounter{figure}\label{cross}\addtocounter{figure}{-1}
\begin{center}
\includegraphics[width=2truein]{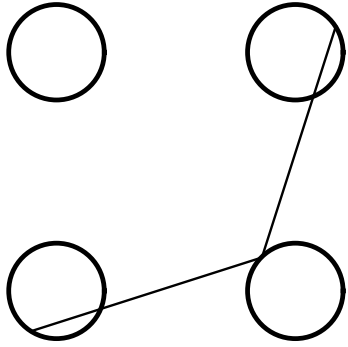}
\caption{A trajectory piece crossing the first and last obstacles}
\end{center}
\end{figure}

\begin{corollary}\label{shortest}
If the trajectory piece from Lemma~\ref{unique} exists and has
admissible type, then it is the
shortest path of type $B$ starting at $\vecx_0$ and ending at $\vecx_s$.
\end{corollary}

\begin{proof}
Similarly as in the proof of Theorem~\ref{thmadm}, the shortest path
of type $B$ from $\vecx_0$ to $\vecx_s$ is a trajectory piece (here we allow
the first segment of the trajectory piece to cross $O_{\veck_0}$ and the
last one to cross $O_{\veck_s}$). By Lemma~\ref{unique}, it is equal to
the trajectory piece from that lemma.
\end{proof}

Of course this trajectory piece depends on $\vecx_0$ and $\vecx_s$. However,
its length depends on those two points only up to an additive
constant. Denote by $c$ the diameter of $O$.

\begin{lemma}\label{const}
For every admissible finite sequence $B=(\veck_n)_{n=0}^s$ of elements of $\Z^m$ the
lengths of trajectory pieces of type $B$ (even if we allow them to
cross $O_{\veck_0}$ and $O_{\veck_s}$) differ by at most $2c$.
The displacements along those trajectory pieces also differ by at most
$2c$.
\end{lemma}

\begin{proof}
Let $\Gamma$ and $\Gamma'$ be two such trajectory pieces, joining
$\vecx_0$ with $\vecx_s$ and $\vecy_0$ with $\vecy_s$ respectively, where $\vecx_0,\vecy_0$
belong to the boundary of $O_{\veck_0}$ and $\vecx_s,\vecy_s$ belong to the
boundary of $O_{\veck_s}$. Replace the first segment of $\Gamma$ by adding
to it the segment joining $\vecx_0$ with $\vecy_0$, and do similarly with the
last segment of $\Gamma$. Then we get a path joining $\vecy_0$ with $\vecy_s$
of type $B$. By Corollary~\ref{shortest}, its length is not smaller
than the length of $\Gamma'$. On the other hand, its length is not
larger than the length of $\Gamma$ plus $2c$. Performing the same
construction with the roles of $\Gamma$ and $\Gamma'$ reversed, we
conclude that the difference of the lengths of those two paths is not
larger than $2c$.

The second statement of the lemma is obvious.
\end{proof}

One can look at the definition of an admissible sequence in the
following way. Instead of a sequence $(\veck_n)_{n=0}^\infty$ of elements
of $\Z^m$ we consider the sequence $(\vecl_n)_{n=1}^\infty$, where
$\vecl_n=\veck_n-\veck_{n-1}$. Since $\veck_0=\zero$, knowing
$(\vecl_n)_{n=1}^\infty$ we can recover $(\veck_n)_{n=0}^\infty$. Now,
condition (\ref{adm3}) can be restated as no obstacle between
$O_\zero$ and $O_{\vecl_n}$, and condition (\ref{adm4}) as the obstacle
$O_{\vecl_n}$ not between $O_\zero$ and $O_{\vecl_n+\vecl_{n+1}}$. Let
$G$ be the directed graph whose vertices are those
$\vecj\in\Z^m\setminus\{\zero\}$ for which there is no obstacle between
$O_\zero$ and $O_\vecj$, and there is an edge (arrow) from $\vecj$ to $\veci$ if
and only if $O_\vecj$ is not between $O_\zero$ and $O_{\vecj+\veci}$. Then every
sequence $(\vecl_n)_{n=1}^\infty$ obtained from an admissible sequence
is a one-sided infinite path in $G$, and vice versa, each one-sided
infinite path in $G$ is a sequence $(\vecl_n)_{n=1}^\infty$ obtained
from an admissible sequence. Hence, we can speak about paths
\emph{corresponding} to admissible sequences and admissible sequences
\emph{corresponding} to paths.

\begin{lemma}\label{fingr}
The set of vertices of $G$ is finite.
\end{lemma}

\begin{proof}
Fix an interior point $\vecx$ of $O_\zero$. By Lemma~\ref{infref}, any ray
beginning at $\vecx$ intersects the interior of some $O_\veck$ with $\veck\ne\zero$.
Let $V_\veck$ be the set of directions (points of the unit sphere) for
which the corresponding ray intersects the interior of $O_\veck$. This set
is open, so we get an open cover of a compact unit sphere. It has a
finite subcover, so there exists a constant $M>0$ such that every ray
from $\vecx$ of length $M$ intersects the interior of some $O_\veck$ with
$\veck\ne\zero$. This proves that the set of vertices of $G$ is finite.
\end{proof}

Note that in $G$ there is never an edge from a vertex to itself.
Moreover, there is a kind of symmetry in $G$. Namely, if $\veck$ is a
vertex then $-\veck$ is a vertex; there is an edge from $\veck$ to $-\veck$; and
if there is an edge from $\veck$ to $\vecj$ then there is an edge from $-\vecj$ to
$-\veck$.

The following lemma establishes another symmetry in $G$.

\begin{lemma}\label{st2}
If $\veck,\vecj\in\Z^m$ and $O_\veck$ is between $O_\zero$ and $O_{\veck+\vecj}$, then
$O_\vecj$ is also between $O_\zero$ and $O_{\veck+\vecj}$. Thus, if there is an
edge in $G$ from $\veck$ to $\vecj$ then there is an edge from $\vecj$ to $\veck$.
\end{lemma}

\begin{proof}
The map $f(x)=\veck+\vecj-\vecx$ defines an isometry of $\Z^m$ and
$f(O_\zero)=O_{\veck+\vecj}$, $f(O_{\veck+\vecj})=O_\zero$, $f(O_\veck)=O_\vecj$. This proves
the first statement of the lemma. The second statement follows from
the first one and from the definition of edges in $G$.
\end{proof}

We will say that the obstacle $O$ is \emph{small} if it is contained
in a closed ball of radius smaller than $\sqrt{2}/4$. To simplify the
notation, in the rest of the paper, whenever the obstacle is small, we
will be using the lifting to $\R^m$ such that the centers of the balls
of radii smaller than $\sqrt{2}/4$ containing the obstacles will be at
the points of $\Z^m$.

Denote by $U$
the set of unit vectors from $\Z^m$ (that is, the ones with one
component $\pm 1$ and the rest of components 0), and by $A_m$ the set
$\{-1,\,0,\,1\}^m\setminus\{\zero\}$ (we use the subscript $m$ by $A$,
since this set will be used sometimes when we consider all dimensions
at once). In particular, $U\subset A_m$.

\begin{lemma}\label{scalar}
Let $O$ be small. If $\veck,\vecl\in\Z^m\setminus\{\zero\}$ and $\langle
\veck,\vecl\rangle\le 0$, then $O_\veck$ is not between $O_\zero$ and
$O_{\veck+\vecl}$. In particular, if $\veck$ and $\vecl$ are vertices of $G$ and
$\langle\veck,\vecl\rangle\le 0$, then there are edges in $G$ from $\veck$ to
$\vecl$ and from $\vecl$ to $\veck$.
\end{lemma}

\begin{proof}
We will use elementary geometry. Consider the triangle with vertices
$A=\zero$, $B=\veck$ and $C=\veck+\vecl$. The angle at the vertex $B$ is at
most $\pi/2$, and the lengths of the sides $AB$ and $BC$ are at least
1. We need to construct a straight line which separates the plane $P$
in which the triangle $ABC$ lies into two half-planes with the first
one containing the open disk of radius $\sqrt{2}/4$ centered at $B$
and the second one containing such disks centered at $A$ and $C$. Then
the hyperplane of dimension $m-1$ through this line and perpendicular
to the plane $P$ will separate $O_\veck$ from $O_\zero$ and $O_{\veck+\vecl}$.
This will prove that there is an edge in $G$ from $\veck$ to $\vecl$. By
Lemma~\ref{st2}, there will be also an edge in $G$ from $\vecl$ to $\veck$.

Let $D$ and $E$ be the points on the sides $BA$ and $BC$ respectively,
whose distance from $B$ is $1/2$ and let $L$ be the straight line
through $D$ and $E$. Since the angle at the vertex $B$ is at most
$\pi/2$, the distance of $B$ from $L$ is at least $\sqrt{2}/4$. Since
$|AD|\ge|BD|$ and $|CE|\ge|BE|$, the distances of $A$ and $C$ from $L$
are at least as large as the distance of $B$ from $L$. This completes
the proof.
\end{proof}

\begin{lemma}\label{am}
For a billiard on a torus with a small obstacle, all elements of $A_m$
are vertices of $G$.
\end{lemma}

\begin{proof}
Let $\vecu\in A_m$ and $\vecv\in\Z^m\setminus\{\zero,\vecu\}$. If
$\vecv=(v_1,v_2,\dots,v_m)$ then $|\langle\vecv,\vecu\rangle|\le\sum_{i=1}^m
|v_i|\le\|\vecv\|^2$, so $\langle\vecv,\vecu-\vecv\rangle\le 0$. Therefore, by
Lemma~\ref{scalar}, $\vecv$ is not between $\zero$ and $u$. This proves
that $u$ is a vertex of $G$.
\end{proof}

\begin{lemma}\label{small}
Assume that $O$ is small. Then $G$ is connected, and for every
vertices $\veck,\vecl$ of $G$ there is a path of length at most $3$ from
$\veck$ to $\vecl$ in $G$, via elements of $U$.
\end{lemma}

\begin{proof}
By Lemma~\ref{am}, the set of vertices of $G$ contains $U$. Let
$\veck,\vecl$ be vertices of $G$. Then $\veck,\vecl\ne\zero$, so there exist
elements $\vecu,\vecv$ of $U$ such that $\langle\veck,\vecu\rangle\le 0$ and
$\langle\vecl,\vecv\rangle\le 0$. By Lemma~\ref{scalar}, there are edges
in $G$ from $\veck$ to $\vecu$ and from $\vecv$ to $\vecl$. If $\vecu=\vecv$ then $\veck\vecu\vecl$
is a path of length 2 from $\veck$ to $\vecl$. If $\vecu\ne\vecv$ then $\langle
\vecu,\vecv\rangle=0$, so by Lemma~\ref{scalar} there is an edge from $\vecu$ to
$\vecv$. Then $\veck\vecu\vecv\vecl$ is a path of length 3 from $\veck$ to $\vecl$.
\end{proof}

\section{Rotation set - torus}\label{sec_rotset}

Now we have enough information in order to start investigating the
\emph{rotation set} $R$ of our billiard. It consists of limits of the sequences
$((\vecy_n-\vecx_n)/t_n)_{n=1}^\infty$, where there is a trajectory piece in
the lifting from
$\vecx_n$ to $\vecy_n$ of length $t_n$, and $t_n$ goes to infinity. Since we
have much larger control of pieces of trajectories of admissible type,
we introduce also the \emph{admissible rotation set} $AR$, where in
the definition we consider only such pieces. Clearly, the admissible
rotation set is contained in the rotation set. By the definition, both
sets are closed. It is also clear that they are contained in the
closed unit ball in $\R^m$, centered at the origin. Due to the
time-reversibility, both sets $R$ and $AR$ are centrally symmetric
with respect to the origin.

For a given point $p$ in the phase space let us consider the
trajectory $t\mapsto T(t)$ in $\R^m$ starting at $p$. We
can ask whether the limit of $(T(t)-T(0))/t$, as $t$ goes to infinity,
exists. If it does, we will call it the \emph{rotation vector} of $p$.
Clearly, it is the same for every point in the phase space of the full
trajectory of $p$, so we can speak of the \emph{rotation vector of a
trajectory}. In particular, every periodic orbit has a rotation vector,
and it is equal to $(T(s)-T(0))/s$, where $s$ is the period of the
orbit.

Note that if we use the discrete time (the number of reflections)
rather than continuous time, we would get all good properties of the
admissible rotation set from the description of the admissible
sequences via the graph $G$ and the results of \cite{Z}. Since we are
using continuous time, the situation is more complicated.
Nevertheless, Lemma~\ref{const} allows us to get similar results. For
a trajectory piece $T$ we will denote by $|T|$ its length and by
$\vecd(T)$ its displacement.

\begin{theorem}\label{convex}
The admissible rotation set of a billiard on a torus with a small
obstacle is convex.
\end{theorem}

\begin{proof}
Fix vectors $\vecu,\vecv\in AR$ and a number $t\in(0,1)$. We want to show that
the vector $t\vecu+(1-t)\vecv$ belongs to $AR$. Fix $\eps>0$. By the
definition, there are finite admissible sequences $A,B$ and trajectory pieces
$T,S$ of type $A,B$ respectively, such that
\begin{equation}\label{eq2}
\left\|\frac{\vecd(T)}{|T|}-\vecu\right\|<\eps \aand \left\|\frac{\vecd(S)}
{|S|}-\vecv\right\|<\eps.
\end{equation}
Both $A,B$ can be represented as finite paths in the graph $G$. By
Lemma~\ref{small}, there are admissible sequences $C_1,C_2,C_3$
represented in $G$ as paths of length at most 3, via elements of $U$,
such that the concatenations of the form
$$D=AC_1AC_1\dots AC_1AC_2BC_3BC_3\dots BC_3B$$
are admissible. There exists a trajectory piece $Q$ of type $D$. We
will estimate its displacement and length.

Assume that in $D$ the block $A$ appears $p$ times and the block $B$
appears $q-p$ times. Let $\vecd_A,\vecd_B$ be the total displacements due to
the blocks $A,B$ respectively. We get
$$\|\vecd_A-p\vecd(T)\|\le 2pc \aand \|\vecd_B-(q-p)\vecd(S)\|\le 2(q-p)c.$$
The displacement due to each of the blocks $C_1,C_2,C_3$ is at most of
norm $2+2c$, so the total displacement due to all those blocks is
at most of norm $q(2+2c)$. If we replace all displacements by the
trajectory lengths, we get the same estimates (we use here
Lemma~\ref{const}). Thus we get the following estimates:
$$\|\vecd(Q)-\vecalpha\|\le 4qc+2q \aand
\big||Q|-\beta\big|\le 4qc+2q,$$
where
$$\vecalpha=p\vecd(T)+(q-p)\vecd(S) \aand \beta=p|T|+(q-p)|S|.$$
Therefore
\begin{eqnarray}\label{eq3}
\left\|\frac{\vecd(Q)}{|Q|}-\frac{\vecalpha}{\beta}\right\|
&\le&\left\|\frac{\vecd(Q)}{|Q|}-\frac{\vecalpha}{|Q|}\right\|
+\left\|\frac{\vecalpha}{|Q|}-\frac{\vecalpha}{\beta}\right\|\nonumber\\
&\le&\frac{4qc+2q}{|Q|}+\|\vecalpha\|\frac{4qc+2q}{|Q|\beta}\\
&=&(4c+2)\frac{q}{|Q|}\left(1+\frac{\|\vecalpha\|}{\beta}\right).
\nonumber
\end{eqnarray}

Set $s=p|T|/\beta$. Then $1-s=(q-p)|S|/\beta$, so
\begin{equation}\label{eq4}
\frac{\vecalpha}{\beta}=s\frac{\vecd(T)}{|T|}+(1-s)\frac{\vecd(S)}{|S|}.
\end{equation}
By (\ref{eq2}), we get
$$\frac{\|\vecalpha\|}{\beta}\le\max(\|\vecu\|,\|\vecv\|)+\eps.$$
Moreover,
$$\frac{|Q|}{q}\ge\frac{\beta}{q}-(4c+2)\ge\min(|T|,|S|)-
(4c+2).$$
Therefore if $|T|$ and $|S|$ are sufficiently large (we may assume
this), the right-hand side of (\ref{eq3}) is less than $\eps$.
Together with (\ref{eq4}), we get
$$\left\|\frac{\vecd(Q)}{|Q|}-\left(s\frac{\vecd(T)}{|T|}+(1-s)\frac{\vecd(S)}
{|S|}\right)\right\|<\eps.$$
By this inequality and (\ref{eq2}), it remains to show that by the
right choice of $p,q$ we can approximate $t$ by $s$ with an arbitrary
accuracy.

We can write $s=f(x)$, where $x=p/q$ and
$$f(x)=\frac{|T|x}{|T|x+|S|(1-x)}.$$
The function $f$ is continuous on $[0,1]$, takes value 0 at 0 and
value 1 at 1. Therefore the image of the set of rational numbers from
$(0,1)$ is dense in $[0,1]$. This completes the proof.
\end{proof}

\begin{theorem}\label{perdense}
For a billiard on a torus with a small obstacle, rotation vectors of
periodic orbits of admissible type are dense in the admissible
rotation set.
\end{theorem}

\begin{proof}
Fix a vector $\vecu\in AR$ and $\eps>0$. We want to find a periodic orbit
of admissible type whose rotation vector is in the $\eps$-neighborhood
of $\vecu$. By the definition, there is an admissible sequence $A$ and a
trajectory piece $T$ of type $A$ such that
\begin{equation}\label{eq5}
\left\|\frac{\vecd(T)}{|T|}-\vecu\right\|<\frac{\eps}{2}.
\end{equation}
Moreover, we can assume that $|T|$ is as large as we need. As in the
proof of Theorem~\ref{convex}, we treat $A$ as a path in the graph $G$
and find an admissible sequence $C$ represented in $G$ as paths of
length at most 3, via elements of $U$, such
that the periodic concatenation $D=ACACAC\dots$ is admissible. There
exists a periodic orbit of type $D$. Let $Q$ be its piece
corresponding to the itinerary $AC$. We will estimate its displacement and length.

Similarly as in the proof of Theorem~\ref{convex}, we get
$$|\vecd(Q)-\vecd(T)|\le 4c+2 \aand \big||Q|-|T|\big|\le 4c+2.$$
Therefore
\begin{eqnarray*}
\left\|\frac{\vecd(Q)}{|Q|}-\frac{\vecd(T)}{|T|}\right\|
&\le&\left\|\frac{\vecd(Q)}{|Q|}-\frac{\vecd(T)}{|Q|}\right\|
+\left\|\frac{\vecd(T)}{|Q|}-\frac{\vecd(T)}{|T|}\right\|\\
&\le&\frac{4c+2}{|T|-(4c+2)}+\frac{\|\vecd(T)\|}{|T|}\cdot
\frac{4c+2}{|T|-(4c+2)}.
\end{eqnarray*}
If $|T|$ is sufficiently large then the right-hand side of this
inequality is smaller than $\eps/2$. Together with (\ref{eq5}) we get
$$\left\|\frac{\vecd(Q)}{|Q|}-\vecu\right\|<\eps.$$
This completes the proof.
\end{proof}

We will refer to closed paths in $G$ as \emph{loops}.

\begin{remark}\label{vertex}
It is clear that in the above theorem we can additionally require that
the corresponding loop in the graph $G$ passes through a given vertex.
\end{remark}

To get more results, we need a generalization of a lemma from
\cite{MZ} to higher dimensions.

\begin{lemma}\label{geom}
Assume that $\zero\in\R^m$ lies in the interior of the convex hull of
a set of $m+1$ vectors $\vecv_0,\vecv_1,\dots,\vecv_m$. For every $\veck>0$
if $L$ is large enough then the following property holds. If $x\in\R^m$ and
$\|x\|\le L$ then there exists $i\in\{0,1,\dots,m\}$ and a positive
integer $n$ such that $\|\vecx+n\vecv_i\|\le L-K$. Moreover, $\|\vecx+j\vecv_i\|\le L$
for $j=1,2,\dots,n-1$.
\end{lemma}

\begin{proof}
Let us fix $\veck>0$. We will consider only $L$ such that $L>K$. Set
$M=\max_i\|\vecv_i\|$.

For each $\vecx\in\R^m$ with $\|\vecx\|=1$ let $f(\vecx)$ be the minimum of
$\|\vecx+t\vecv_i\|$ over $i=0,1,\dots,m$ and $t\ge 0$. By the assumption,
$f(\vecx)<1$. Clearly $f$ is continuous, and therefore there is $\eps>0$
such that $f(\vecx)\le 1-\eps$ for every $\vecx$. Thus, for every $\vecy\in\R^m$
there is $i=0,1,\dots,m$ and $s\ge 0$ such that
$\|\vecy+s\vecv_i\|\le(1-\eps)\|\vecy\|$. Let $n$ be the smallest integer larger
than $s$. Then $n>0$, and if $L\ge(M+K)/\eps$ and $\|\vecy\|\le L$ then
$\|\vecy+n\vecv_i\|\le(1-\eps)L+M\le L-K$.

The last statement of the lemma follows from the convexity of the
balls in $\R^m$.
\end{proof}

Now we can follow the methods of \cite{MZ} and \cite{Z}. We assume
that our billiard has a small obstacle. For a full trajectory $T$ we will
denote by $T(t)$ the point to which we get after time $t$.

\begin{lemma}\label{follow}
If $\vecu$ is a vector from the interior of $AR$, then there exists a
full trajectory $T$ of admissible type and a constant $M$ such that
\begin{equation}\label{eq6}
\|T(t)-T(0)-t\vecu\|\le M
\end{equation}
for all $t\in\R$.
\end{lemma}

\begin{proof}
Let us think first of positive $t$'s.
Since $\vecu$ is in the interior of $AR$, one can choose $m+1$ vectors
$\vecw_0,\vecw_1,\dots,\vecw_m\in AR$ such that $u$ is in the interior of the
convex hull of those vectors. Moreover, by Theorem~\ref{perdense} and
Remark~\ref{vertex} we may assume that $\vecw_i$ are rotation vectors of
periodic orbits $P_i$ of admissible type, corresponding to loops $A_i$ in $G$
passing through a common vertex $V$. We can also consider those loops
as finite paths, ending at $V$ and starting at the next vertex in the
loop. Set
$$\vecv_i=\vecd(P_i)-|P_i|\vecu=|P_i|\vecw_i-|P_i|\vecu=|P_i|(\vecw_i-\vecu).$$
Since $\vecu$ is in the interior of the convex hull of the vectors $\vecw_i$,
we get that $\zero$ is in the interior of the convex hull of the
vectors $\vecw_i-\vecu$, and therefore $\zero$ is in the interior of the
convex hull of the vectors $\vecv_i$.

We will construct our trajectory, or rather a corresponding path in
the graph $G$, by induction, using Lemma~\ref{geom}. Then we get a
corresponding trajectory of admissible type by Theorem~\ref{thmadm}. We start with
the empty sequence, that corresponds to the trajectory piece
consisting of one point. Then, when a path $B_j$ in $G$
(corresponding to a trajectory piece $Q_j$) is constructed, and it
ends at $V$, we look at the vector $\vecx=\vecd(Q_j)-|Q_j|\vecu$ and choose $\vecv_i$
and $n$ according to Lemma~\ref{geom}. We append $B_j$ by adding $n$
repetitions of $A_i$ (corresponding to a trajectory piece that we can
call $nP_i$) and obtain $B_{j+1}$(corresponding to a trajectory piece
$Q_{j+1}$). To do all this, we have to define $\veck$ that is used in
Lemma~\ref{geom} and prove that if $\|\vecx\|\le L$ then also
$\big\|\vecd(Q_{j+1})-|Q_{j+1}|\vecu\big\|\le L$.

Let us analyze the situation. When we concatenate $Q_j$ and
$A_i\dots A_i$ ($n$ times) to get $Q_{j+1}$, by Lemma~\ref{const}
we have
$$\|\big(\vecd(Q_{j+1})-|Q_{j+1}|\vecu\big)-\big(\vecd(Q_j)-|Q_j|\vecu\big)-\big(\vecd(nP_i)-|nP_i|\vecu\big)\|\le
4c(1+\|\vecu\|).$$
Moreover,
$$\vecd(nP_i)-|nP_i|\vecu=n(\vecd(P_i)-|P_i|\vecu)=n\vecv_i.$$
Therefore in Lemma~\ref{geom} we have to take $\veck=4c(1+\|\vecu\|)$ and then we can
make the induction step.
In such a way we obtain an infinite path $B$ in $G$. By
Theorem~\ref{thmadm}, there exists a billiard trajectory $T$ of type
$B$.

Note that we did not complete the proof yet, because we got
(\ref{eq6}) (with $M=L$) only for a sequence of times $t=|Q_j|$. We
can do better using the last statement of Lemma~\ref{geom}. This shows that
(\ref{eq6}) with $M=L+K$ holds for a sequence of times $t$ with the
difference of two consecutive terms of this sequence not exceeding
$s=\max(|P_0|,|P_1|,\dots,|P_n|)+4c$. Every time $t'$ can be written
as $t+r$ with $t$ being a term of the above sequence (so that
(\ref{eq6}) holds with with $M=L+K$) and $r\in[0,s)$. Then
$$\|T(t')-T(0)-t'\vecu\|\le L+K+\|T(t+r)-T(t)\|+r\|\vecu\|.$$
Thus, (\ref{eq6}) holds for all times with $M=L+K+s+s\|\vecu\|$.

The same can be done for negative $t$'s, so we get a full (two-sided)
path, and consequently a full trajectory.
\end{proof}

Now we are ready to prove the next important theorem. Remember that
our phase space is a factor of a compact connected subset of the unit
tangent bundle over the torus.

\begin{theorem}\label{cominv}
For a billiard on a torus with a small obstacle, if $\vecu$ is a vector
from the interior of $AR$, then there exists a compact
invariant subset $Y$ of the phase space, such that every trajectory
from $Y$ has admissible type and rotation vector $\vecu$.
\end{theorem}

\begin{proof}
Let $Y$ be the closure of the trajectory $T$ from Lemma~\ref{follow},
taken in the phase space. If $S$ is a trajectory obtained from $T$ by
starting it at time $s$ (that is, $S(t)=T(s+t)$) then by
Lemma~\ref{follow} we get $\|S(t)-S(0)-t\vecu\|\le 2M$ for all $t$.
By continuity of the flow, this property extends to every trajectory
$S$ from $Y$. This proves that every trajectory from $Y$ has rotation
vector $\vecu$.

Since a trajectory of admissible type has no tangencies to the
obstacle (by the condition~\ref{adm4} of the definition of admissible
sequences), so each finite piece of a trajectory from $Y$ has
admissible type. Therefore every trajectory from $Y$ has admissible
type.
\end{proof}

\begin{remark}\label{minimal}
The set $Y$ above can be chosen minimal, and therefore the trajectory
from Lemma~\ref{follow} can be chosen recurrent.
\end{remark}

As a trivial corollary to Theorem~\ref{cominv} we get the following.

\begin{corollary}\label{pointwise}
For a billiard on a torus with a small obstacle, if $\vecu$ is a vector
from the interior of $AR$, then there exists a trajectory of
admissible type with rotation vector $\vecu$.
\end{corollary}

We also get another corollary, which follows from the existence of an
ergodic measure on $Y$.

\begin{corollary}\label{ergodic}
For a billiard on a torus with a small obstacle, if $\vecu$ is a vector
from the interior of $AR$, then there exists an ergodic invariant
probability measure in the phase space, for which the integral of the
velocity is equal to $\vecu$ and almost every trajectory is of
admissible type.
\end{corollary}

This corollary is stronger than Corollary~\ref{pointwise}, because
from it and from the Ergodic Theorem it follows that almost every
point has rotation vector $\vecu$. The details of the necessary formalism
are described in Section~\ref{sec_conn}. Of course, in our particular
case both results are corollaries to Theorem~\ref{cominv}, so we know
anyway that all points of $Y$ have rotation vector $\vecu$.

\section{Admissible rotation set is large}\label{sec_arlarge}

In this section we will investigate how large the admissible rotation
set $AR$ is. This of course depends on the size of the obstacle and
the dimension of the space. We will measure the size of $AR$ by the
radius of the largest ball centered at the origin and contained in $AR$.

We will start with the estimates that depend on the dimension $m$ of
the space but not on the size of the obstacle (provided it is small in
our meaning). In order to do it, we first identify some elements of
$\Z^m$ that are always vertices of $G$. Set
$$A_m=\{-1,\,0,\,1\}^m\setminus\{\zero\}.$$

\begin{lemma}\label{amar}
If $\veck\in A_m$ then $(\sqrt{2}/2)(\veck/\|\veck\|)\in AR$.
\end{lemma}

\begin{proof}
If $\veck\in U$ then there is a vector $\vecl\in U$ orthogonal to $\veck$.
Vectors $\veck+\vecl$ and $\veck-\vecl$ belong to $A_m$ and one can easily check
that there are edges from $\veck+\vecl$ to $\veck-\vecl$ and from $\veck-\vecl$ to
$\veck+\vecl$ in $G$. The periodic path $(\veck+\vecl)(\veck-\vecl)(\veck+\vecl)(\veck-\vecl)
\dots$ in $G$ gives us a periodic orbit $P$ of the billiard. The
displacement along $P$ is $2\veck$ and the period of $P$ is smaller than
$\|\veck+\vecl\|+\|\veck-\vecl\|=2\sqrt{2}$, so the rotation vector of $P$ is
$t\veck$, where $t>\sqrt{2}/2$. Since $\zero\in AR$ and $AR$ is convex, we
get $(\sqrt{2}/2)(\veck/\|\veck\|)\in AR$.

Assume now that $\veck\in A_m$ and $\|\veck\|>1$. Then $\veck=\vecl+\vecu$ for some
$\vecl\in A_m$ and $\vecu\in U$ such that $\vecu$ is orthogonal to $\vecl$. By
Lemma~\ref{am}, $\vecl$ is a vertex of $G$. By Lemma~\ref{scalar}
there are edges in $G$ from $\vecl$ to $\vecu$ and from $\vecu$ to $\vecl$.
Similarly as before, we get a periodic orbit of the billiard
(corresponding to the periodic path $\vecl\vecu\vecl\vecu\dots$) with the
displacement $\veck$ and period less than $\|\vecl\|+\|\vecu\|=\sqrt{\|\veck\|^2-1}
+1$, so
$$\frac{\|\veck\|}{\sqrt{\|\veck\|^2-1}+1}\cdot\frac{\veck}{\|\veck\|}\in AR.$$
Since $\|\veck\|/(\sqrt{\|\veck\|^2-1}+1)\ge\sqrt{2}/2$, the vector
$(\sqrt{2}/2)(\veck/\|\veck\|)$ also belongs to $AR$.
\end{proof}

By the results of \cite{Nandor}, the convex hull of $A_m$ contains the
closed ball centered at $\zero$ with radius $2/\sqrt{\ln m+5}$. From
this and Lemma~\ref{amar} we get immediately the following result.

\begin{theorem}\label{large}
For a billiard on a torus with a small obstacle, the set $AR$ contains
the closed ball centered at $\zero$ with radius $\sqrt{2/(\ln m+5)}$.
\end{theorem}

Now we proceed to the estimates that are independent of the dimension
$m$. This is not as simple as it seems. As we saw above, a
straightforward attempt that takes into account only those vectors of
$\Z^m$ for which we can show explicitly that they are vertices of $G$,
gives estimates that go to $0$ as $m\to\infty$. By the results of
\cite{Nandor}, those estimates cannot be significantly improved.
Therefore we have to use another method.

Let us assume first that $O_\zero$ is the ball centered at $\zero$ of
radius $r<\sqrt{2}/4$. We start with a simple lemma.

\begin{lemma}\label{st1}
Assume that $O_\vecl$ is between $O_\zero$ and $O_\veck$ and let $\theta$
be the angle between the vectors $\veck$ and $\vecl$. Then
\begin{equation}\label{eqst1}
\langle\veck,\vecl\rangle^2\ge\|\veck\|^2\left(\|\vecl\|^2-4r^2\right)>
\|\veck\|^2\left(\|\vecl\|^2-\frac{1}{2}\right).
\end{equation}
and
\begin{equation}\label{eqst2}
\sin\theta\le 2r/\|\vecl\|.
\end{equation}
\end{lemma}

\begin{proof}
If $O_\vecl$ is between $O_\zero$ and $O_\veck$ then there is a line
parallel to the vector $\veck$, whose distances from $\zero$ and $\vecl$
are at most $r$. Therefore the distance of $\vecl$ from the line
through $\zero$ and $\veck$ is at most $2r$. The orthogonal projection
of $\vecl$ to this line is $(\langle\veck,\vecl\rangle/\|\veck\|^2)k$, so
$$\left\|\vecl-\frac{\langle\veck,\vecl\rangle}{\|\veck\|^2}\veck\right\|^2\le 4r^2.$$
The left-hand side of this inequality is equal to
$$\|\vecl\|^2-\frac{\langle\veck,\vecl\rangle^2}{\|\veck\|^2}$$
and $4r^2<1/2$, so (\ref{eqst1}) holds.

By (\ref{eqst1}), we have
$$\sin^2\theta=1-\cos^2\theta=1-\frac{\langle\veck,\vecl\rangle^2}
{\|\veck\|^2\|\vecl\|^2}\le 1-\frac{\|\veck\|^2(\|\vecl\|^2-4r^2)}
{\|\veck\|^2\|\vecl\|^2}=\frac{4r^2}{\|\vecl\|^2},$$
so (\ref{eqst2}) holds.
\end{proof}

Clearly, the angle $\theta$ above is acute.

The estimate in the next lemma requires extensive use of the fact that
the vectors that we are considering have integer components.

\begin{lemma}\label{st3}
Assume that $O_\vecl$ is between $O_\zero$ and $O_\veck$ and
$$\langle\veck,\vecl\rangle\le\langle\veck,\veck-\vecl\rangle.$$
Then $\|\vecl\|\le\|\veck\|/2$.
\end{lemma}

\begin{proof}
By Lemma~\ref{st1}, (\ref{eqst1}) holds. Since $\langle\veck,\vecl\rangle\le
\langle\veck,\veck-\vecl\rangle$, we get $\langle\veck,\vecl\rangle\le\|\veck\|^2-
\langle\veck,\vecl\rangle$, so
\begin{equation}\label{eqst3}
2\langle\veck,\vecl\rangle\le\|\veck\|^2.
\end{equation}
By (\ref{eqst1}) and (\ref{eqst3}) we get
\begin{equation}\label{eqst4}
\|\veck\|^2>4\|\vecl\|^2-2.
\end{equation}
If $\|\vecl\|>\|\veck\|/2$, then $\|\veck\|^2<4\|\vecl\|^2$. Together with
(\ref{eqst4}), since $\|\veck\|^2$ and $\|\vecl\|^2$ are integers, we get
\begin{equation}\label{eqst5}
\|\veck\|^2=4\|\vecl\|^2-1.
\end{equation}
 From (\ref{eqst3}) and (\ref{eqst5}) we get $\langle\veck,\vecl\rangle\le
2\|\vecl\|^2-1/2$. Hence, since $\langle\veck,\vecl\rangle$ is also an
integer, we get $\langle\veck,\vecl\rangle\le 2\|\vecl\|^2-1$. From this,
(\ref{eqst1}) and (\ref{eqst5}) we get
$$(2\|\vecl\|^2-1)^2>(4\|\vecl\|^2-1)\left(\|\vecl\|^2-\frac{1}{2}\right)
=\left(2\|\vecl\|^2-\frac{1}{2}\right)(2\|\vecl\|^2-1),$$
a contradiction.
\end{proof}

Let us think about standing at the origin and looking at the sky,
where vertices of $G$ are stars. Are there big parts of the sky
without a single star? Observe that as the dimension $m$
of the space grows, the angles between the integer vectors tend to
become larger. For instance, the angle between the vectors
$(1,0,\dots,0)$ and $(1,\dots,1)$ is of order $\pi/2-1/\sqrt{m}$.
Thus, any given acute angle can be considered relatively small if $m$
is sufficiently large.

Let $\alpha$ be a positive angle. We will say that a set
$A\subset\Z^m\setminus\{\zero\}$ is \emph{$\alpha$-dense in the sky}
if for every $\vecv\in\R^m\setminus\{\zero\}$ there is $\vecu\in A$ such that
the angle between the vectors $\vecv$ and $\vecu$ is at most $\alpha$. Set
$$\eta(r)=\sum_{n=0}^\infty\arcsin\frac{r}{2^{n-1}}.$$

\begin{proposition}\label{st5}
The set of vertices of $G$ is $\eta(r)$-dense in the sky.
\end{proposition}

\begin{proof}
Fix a vector $\vecv\in\R^m\setminus\{\zero\}$ and $\eps>0$.
There exists $\veck_0\in\Z^m\setminus\{\zero\}$ such that
the angle between $\vecv$ and $\veck_0$ is less than $\eps$. Then we define by
induction a finite sequence $(\veck_1,\veck_2,\dots,\veck_n)$ of elements of
$\Z^m\setminus\{\zero\}$ such that $O_{\veck_{i+1}}$ is between $O_\zero$
and $O_{\veck_i}$ and $\|\veck_{i+1}\|\le\|\veck_i\|/2$. This is possible by
Lemmas~\ref{st2} and~\ref{st3}. Since $\|\veck_i\|\ge 1$ for each $i$,
this procedure has to terminate at some $\veck_n$. Then there is no
obstacle between $O_\zero$ and $O_{\veck_n}$, so $\veck_n$ is a vertex of $G$.
By Lemma~\ref{st1}, the angle between $\veck_i$ and $\veck_{i+1}$ is at most
$\arcsin(2r/\|\veck_{i+1}\|)$. By our construction, we have $\|\veck_n\|\ge
1=2^0$, $\|\veck_{n-1}\|\ge 2^1$, $\|\veck_{n-2}\|\ge 2^2$, etc. Therefore the
angle between $\veck_n$ and $\veck_0$ is smaller than $\eta(r)$. Hence, the
angle between $\vecv$ and $\veck_n$ is smaller than $\eta(r)+\eps$. Since
$\eps$ was arbitrary, this angle is at most $\eta(r)$.
\end{proof}

\begin{remark}\label{st6}
Proposition~\ref{st5} was proved under the assumption that $O_\zero$
is the ball centered at $\zero$ of radius $r<\sqrt{2}/4$. However,
making an obstacle smaller results in preservation or even enlargement
of $G$. Moreover, we have a freedom in the lifting where to put the
origin. Therefore, Proposition~\ref{st5} remains true under a weaker
assumption, that $O$ is contained in a closed ball of radius
$r<\sqrt{2}/4$.
\end{remark}

Let us investigate the properties of $\eta(r)$.

\begin{lemma}\label{st7}
The function $\eta$ is continuous and increasing on $(0,\sqrt{2}/4]$.
Moreover,\break $\eta(r)<\sqrt{2}\,\,\pi r$. In particular,
$\eta(\sqrt{2}/4)<\pi/2$ and
$$\lim_{r\to 0}\eta(r)=0.$$
\end{lemma}

\begin{proof}
Assume that $0<r\le\sqrt{2}/4$. Then all numbers whose arcus sine we
are taking are from the interval $(0,\sqrt{2}/2]$, so clearly $\eta$
is continuous and increasing. We have also the estimate
$$\frac{x}{\arcsin x}\ge\frac{\sqrt{2}/2}{\pi/4}=
\frac{2\sqrt{2}}{\pi}.$$
Moreover, the equality holds only if $x=\sqrt{2}/2$. Thus
$$\eta(r)<\sum_{n=0}^\infty\frac{\pi r}{2^n\,\sqrt{2}}=\sqrt{2}\,\,\pi r.$$
Therefore $\lim_{r\to 0}\eta(r)=0$ and
$$\eta(\sqrt{2}/4)<\sqrt{2}\,\,\pi\frac{\sqrt{2}}{4}=\frac{\pi}{2}.$$
\end{proof}

Now we assume only that $O$ is small. In the next lemma we obtain two
estimates of the length of $t\veck\in AR$ if $\veck$ is a vertex of $G$. One of
those estimates will be useful for all vertices of $G$, the other one
for those with large norm. The main idea of the proof is similar as in
the proof of Lemma~\ref{amar}.

\begin{lemma}\label{st8}
If $\veck$ is a vertex of $G$ then the vectors $(1-\sqrt{2}/2)(\veck/\|\veck\|)$
and $\big((\|\veck\|-1)/(\|\veck\|+1)\big)(\veck/\|\veck\|)$ belong to $AR$.
\end{lemma}

\begin{proof}
Let $\veck=(x_1,x_2,\dots,x_m)$ be a vertex of $G$. Let $s$ be the number
of non-zero components of $\veck$. Then we may assume that $x_i\ne 0$ if
$i\le s$ and $x_i=0$ if $i>s$. If $s=1$ then the statement of the
lemma follows from Lemma~\ref{amar}.

Assume now that $s>1$. Then for every $i\le s$ there is a vector
$\vecv_i\in U$ with only $i$-th component non-zero and $\langle
\vecv_i,\veck\rangle<0$. By Lemma~\ref{scalar} there are edges in $G$ from $\veck$
to $\vecv_i$ and from $\vecv_i$ to $\veck$, so the periodic path $\veck\vecv_i\veck\vecv_i\dots$
in $G$ gives us a periodic orbit $P_i$ of the billiard. The
displacement along $P_i$ is $\veck+\vecv_i$ and the period of $P_i$ is smaller
than $\|\veck\|+\|\vecv_i\|=\|\veck\|+1$, so the rotation vector of $P_i$ is
$t_i(\veck+\vecv_i)$ with $t_i>1/(\|\veck\|+1)$. Therefore the vector
$(\veck+\vecv_i)/(\|\veck\|+1)$ belongs to $AR$.

Since the vectors $\vecv_i$ form an orthonormal basis of $\R^s$, we have
$$\veck=\sum_{i=1}^s\langle \vecv_i,\veck\rangle \vecv_i.$$
Set
$$a=\sum_{i=1}^s\langle \vecv_i,\veck\rangle \aand a_i=\frac{\langle
\vecv_i,\veck\rangle}{a}.$$
Then the vector
$$\vecu=\sum_{i=1}^s a_i\frac{\veck+\vecv_i}{\|\veck\|+1}$$
is a convex combination of elements of $AR$, so $\vecu\in AR$. We have
$$\vecu=\frac{\veck}{\|\veck\|+1}\sum_{i=1}^s a_i+\frac{1}{a(\|\veck\|+1)}\sum_{i=1}^s
\langle \vecv_i,\veck\rangle \vecv_i=\frac{\veck}{\|\veck\|+1}\left(1+\frac{1}{a}\right).$$
For each $i$ we have $\langle \vecv_i,\veck\rangle\le -1$, so $a\le -s$, and
therefore $1+1/a\ge(s-1)/s$. Moreover,
$$\frac{\|\veck\|}{\|\veck\|+1}\ge\frac{\sqrt{s}}{\sqrt{s}+1}.$$
Since $s\ge 2$, we get
$$\frac{s-1}{s}\cdot\frac{\sqrt{s}}{\sqrt{s}+1}=\frac{\sqrt{s}-1}
{\sqrt{s}}\ge\frac{\sqrt{2}-1}{\sqrt{2}}=1-\frac{\sqrt{2}}{2},$$
so the vector $u$ has the direction of $\veck$ and length at least
$1-\sqrt{2}/2$.

To get the other estimate of the length of $\vecu$, note that $\langle
\vecv_i,\veck\rangle=-|x_i|$, so
$$a=-\sum_{i=1}^s|x_i|\le-\|\veck\|,$$
and hence
$$\|\vecu\|\ge\frac{\|\veck\|}{\|\veck\|+1}\cdot\left(1-\frac{1}{\|\veck\|}\right)
=\frac{\|\veck\|-1}{\|\veck\|+1}.$$
\end{proof}

\begin{lemma}\label{st9}
Let $A\subset\R^m$ be a finite set, $\alpha$-dense in the sky for
some $\alpha<\pi/2$. Assume that every vector of $A$ has norm $c$.
Then the convex hull of $A$ contains a ball of radius $c\cos\alpha$,
centered at $\zero$.
\end{lemma}

\begin{proof}
Let $\veck$ be the convex hull of $A$. Then $\veck$ is a convex polytope with
vertices from $A$. Since $A$ is $\alpha$-dense in the sky and
$\alpha<\pi/2$, $\veck$ is non-degenerate and $\zero$ belongs to its
interior. Let $s$ be the radius of the largest ball centered at
$\zero$ and contained in $\veck$. This ball is tangent to some face of $\veck$
at a point $\vecv$. Then the whole $\veck$ is contained in the half-space
$\{\vecu:\langle \vecu,\vecv\rangle\le\|\vecv\|^2\}$. In particular, for every $\vecu\in
A$ we have $\langle \vecu,\vecv\rangle\le\|\vecv\|^2$. Since $A$ is $\alpha$-dense
in the sky, there is $\vecu\in A$ such that the angle between $\vecv$ and $\vecu$
is at most $\alpha$. Therefore
$$\|\vecv\|^2\ge\langle \vecu,\vecv\rangle\ge\|\vecu\|\|\vecv\|\cos\alpha=
c\|\vecv\|\cos\alpha,$$
so $s=\|\vecv\|\ge c\cos\alpha$.
\end{proof}

Now we can get the first, explicit, estimate of the radius of the
largest ball contained in $AR$.

\begin{theorem}\label{st10}
For a billiard on a torus with small obstacle, assume that $O$ is
contained in a closed ball of radius
$r<\sqrt{2}/4$. Then the admissible rotation set contains the closed ball of radius
$(1-\sqrt{2}/2)\cos\eta(r)$ centered at $\zero$.
\end{theorem}

\begin{proof}
Set $H=\{(1-\sqrt{2}/2)(\veck/\|\veck\|):\veck$ is a vertex of $G\}$ and let $\veck$ be the convex
hull of $H$. By Lemma~\ref{st8} and by the convexity of $AR$, we have
$\veck\subset AR$. By Proposition~\ref{st5}, $H$ is $\eta(r)$-dense in the
sky. Thus, by Lemma~\ref{st9}, $\veck$ contains the closed ball of radius
$(1-\sqrt{2}/2)\cos\eta(r)$ centered at $\zero$.
\end{proof}

The second estimate is better for small $r$ (uniformly in $m$), but does not give an
explicit formula for the radius of the ball contained in $AR$. We
first need two lemmas.

\begin{lemma}\label{st11}
For every integer $N>1$ there exists an angle $\beta(N)>0$
(independent of $m$), such that if $\vecu,\vecv\in\Z^m\setminus\zero$ are
vectors of norm less than $N$ and the angle $\theta$ between $\vecu$ and
$\vecv$ is positive, then $\theta\ge\beta(N)$.
\end{lemma}

\begin{proof}
Under our assumptions, each of the vectors $\vecu,\vecv$ has less than $N^2$
non-zero components. Therefore the angle between $\vecu$ and $\vecv$ is the
same as the angle between some vectors $\vecu',\vecv'\in\Z^{2N^2-2}$ of the same
norms as $\vecu,\vecv$. However, there are only finitely many vectors in
$\Z^{2N^2-2}\setminus\{\zero\}$ of norm less than $N$, so the lemma holds
with $\beta(N)$ equal to the smallest positive angle between such
vectors.
\end{proof}

\begin{lemma}\label{st12}
Let $A\subset\Z^m\setminus\zero$ be a finite set, $\alpha$-dense in the sky for some
$\alpha<\beta(N)/2$, where $\beta(N)$ is as in Lemma~\ref{st11}.
Then the set of those elements of $A$ which have norm at least $N$ is
$2\alpha$-dense in the sky.
\end{lemma}

\begin{proof}
Let $B$ be the set of vectors of $\R^m\setminus\{\zero\}$ whose
angular distance from some vector of $A$ of norm at least $N$ is
$2\alpha$ or less.

Let $\vecv\in\R^m$ be a non-zero vector. Assume that the angles between
$\vecv$ and all vectors of $A$ are non-zero. By the assumptions, there
exists $\vecu\in A$ such that the angle $\theta$ between $\vecv$ and $\vecu$ is at
most $\alpha$. If $\|\vecu\|\ge N$ then $\vecv\in B$. Suppose that $\|\vecu\|<N$.
Choose $\eps>0$ such that $\eps<\beta(N)-2\alpha$ and $\eps<\theta$.
We draw a great circle in the sky through $\vecu$ and $\vecv$ and go along it
from $\vecu$ through $\vecv$ and beyond it to some $\vecv'$ so that the angle
between $\vecu$ and $\vecv'$ is $\alpha+\eps$. Now, there exists $\vecu'\in A$
such that the angle between $\vecv'$ and $\vecu'$ is at most $\alpha$. Then
the angle between $\vecu$ and $\vecu'$ is at least $\eps$ and at most
$2\alpha+\eps<\beta(N)$, so $\|\vecu'\|\ge N$. The angle between $\vecv$ and
$\vecu'$ is at most $2\alpha+\eps-\theta<2\alpha$, so again, $\vecv\in B$.

In such a way we have shown that $B$ is dense in
$\R^m\setminus\{\zero\}$. It is also clearly closed in
$\R^m\setminus\{\zero\}$, so it is equal to $\R^m\setminus\{\zero\}$.
\end{proof}

Let $\beta(N)$ be as in Lemma~\ref{st11}, and let $N(r)$ be the
maximal $N$ such that $\eta(r)<\beta(N)/2$. This definition is
correct for sufficiently small $r$, since clearly $\beta(N)\to 0$ as
$N\to\infty$, and by Lemma~\ref{st7}, $\eta(r)\to 0$ as $r\to 0$. It
follows that $N(r)\to\infty$ as $r\to 0$.

\begin{theorem}\label{st13}
For a billiard on a torus with a small obstacle, assume that $O$ is
contained in a closed ball of radius $r<\sqrt{2}/4$
for $r$ so small that $N(r)$ is defined. Then the admissible rotation set contains
the closed ball of radius $\big((N(r)-1)/(N(r)+1)\big)\cos(2\eta(r))$
centered at $\zero$.
\end{theorem}

\begin{proof}
Since $\eta(r)<\beta(N(r))/2$, by
Proposition~\ref{st5}, Remark~\ref{st6} and Lemma~\ref{st12}, the set
of those vertices of $G$ that have norm at least $N(r)$ is
$2\eta(r)$-dense in the sky. By Lemma~\ref{st8}, if $\veck$ is such a
vertex, $\big((\|\veck\|-1)/(\|\veck\|+1)\big)(\veck/\|\veck\|)\in AR$. Since
$\zero\in AR$ and $AR$ is convex, also $\big((N(r)-1)/(N(r)+1)\big)(\veck/\|\veck\|)
\in AR$. Then by Lemma~\ref{st9}, the closed ball of radius
$\big((N(r)-1)/(N(r)+1)\big)\cos(2\eta(r))$ centered at $\zero$ is
contained in $AR$.
\end{proof}

\begin{corollary}\label{st14}
The radius of the largest ball centered at $\zero$ contained in the
admissible rotation set
goes to $1$ uniformly in $m$ as the diameter of the obstacle goes to
$0$.
\end{corollary}

We conclude this section with a result showing that even though $AR$
may be large, it is still smaller than $R$.

\begin{theorem}\label{st15}
For a billiard on a torus with a small obstacle, the admissible
rotation set is contained in the open unit ball. In particular, $AR\ne
R$.
\end{theorem}

\begin{proof}
Since the graph $G$ is finite, there exist positive constants
$c_1<c_2$ such that for every trajectory piece of admissible type the
distance between two consecutive reflections is contained in
$[c_1,c_2]$. Moreover, from the definition of an edge in $G$ and from
the compactness of the obstacle it follows that there is an angle
$\alpha>0$ such that the direction of a trajectory piece of admissible
type changes by at least $\alpha$ at each reflection. Consider the
triangle with two sides of length $c_1$ and $c_2$ and the angle
$\pi-\alpha$ between them. Let $a$ be the ratio between the length of
the third side and $c_1+c_2$, that is,
$$a=\frac{\sqrt{c_1^2+c_2^2+2c_1c_2\cos\alpha}}{c_1+c_2}.$$
This ratio is less than 1 and it decreases when $\alpha$ or $c_1/c_2$
grows. Therefore, if consecutive reflections for a trajectory piece of
admissible type are at times $t_1,t_2,t_3$, then the displacement
between the first and the third reflections divided by $t_3-t_1$ is at
most $a$. Thus, every vector from $AR$ has length at most $a$.

Clearly, the vector $(1,0,\dots,0)$ belongs to $R$, and thus $AR\ne
R$.
\end{proof}

\section{Billiard in the square}\label{sec_square}

Now we consider a billiard in the square $S=[-1/2,1/2]^2$ with one
convex obstacle $O$ with a smooth boundary. The lifting to $\R^2$,
considered in Section~\ref{sec_prelim} is replaced in this case by the
unfolding to $\R^2$. That is, we cover $\R^2$ by the copies of $S$
obtained by consecutive symmetries with respect to the lines $x=n+1/2$
and $y=n+1/2$, $n\in\Z$. Thus, the square $S_\veck=S+\veck$ ($\veck\in\Z^2$) with
the obstacle $O_\veck$ in it is the square $S$ with $O$, translated by
$\veck$, with perhaps an additional symmetry applied. If $\veck=(p,q)$, then,
if both $p,q$ are even, there is no additional symmetry; if $p$ is
even and $q$ odd, we apply symmetry with respect to the line $y=q$; if
$p$ is odd and $q$ is even, we apply symmetry with respect to the line
$x=p$; and if both $p,q$ are odd, we apply central symmetry with
respect to the point $(p,q)$. In this model, trajectories in
$S$ with obstacle $O$ unfold to trajectories in $\R^2$ with
obstacles $O_\veck$, $\veck\in\Z^2$.

The situation in $\R^2$ is now the same as in the case of the torus
billiard, except that, as we mentioned above, the obstacles are not
necessarily the translations of $O_\zero$, and, of course, the
observable whose averages we take to get the rotation set is
completely different. Let us trace which definitions, results and
proofs of Sections~\ref{sec_prelim}, \ref{sec_rotset} and
\ref{sec_arlarge} remain the same, and which need modifications.

The definitions of \emph{between} and \emph{type} remain the same.
Lemma~\ref{infref} is still valid, but in its proof we have to look at
the trajectory on the torus $\R^2/(2\Z)^2$ rather than $\R^2/\Z^2$.
Then the definition of an \emph{admissible sequence} remains the same.

The first part of Theorem~\ref{thmadm} and its proof remains the same,
but in the proof of the part about periodic trajectories we have to be
careful. The point $\vecx_0+\vecp$ from the last paragraph of the proof has to
be replaced by a point that after folding (the operation reverse to the
unfolding) becomes $\vecx_0$. This gives us a periodic orbit in the
unfolding that projects (folds) to a periodic orbit in the square.
Moreover, there may be a slight difference between the square case and
the torus case if we want to determine the least discrete period of
this orbit (where in the square case, in analogy to the torus case, we
count only reflections from the obstacle). In the torus case it is
clearly the same as the least period of the type. In the square case
this is not necessarily so. For instance, if the obstacle is  a disk
centered at the origin, the orbit that goes vertically from the
highest point of the disk, reflects from the upper side of the square
and returns to the highest point of the disk, has discrete period 1 in
the above sense. However, its type is periodic of period 2 and in the
unfolding it has period 2. Fortunately, such things are irrelevant for
the rest of our results.

Lemma~\ref{unique}, Corollary~\ref{shortest} and their proofs remain
the same as in the torus case. The same can be said about the part of
Lemma~\ref{const} that refers to the lengths of trajectory pieces.

The definition of the graph $G$ has to be modified. This is due to the
fact that the conditions~(\ref{adm3}) and~(\ref{adm4}) cannot be restated in
the same way as in the torus case, because now not only translations,
but also symmetries are involved (the obstacle needs not be symmetric,
and the unfolding process involves symmetries about vertical and
horizontal lines). In order to eliminate symmetries, we
enlarge the number of vertices of $G$ four times. Instead of
$O_\zero$, we look at $\{O_\zero,O_{(1,0)},O_{(0,1)},O_{(1,1)}\}$. For
every $\veck=(p,q)\in\Z^2$ there is $\zeta(\veck)\in Q$, where
$Q=\{(0,0),(1,0),(0,1),(1,1)\}$, such that $\veck-\zeta(\veck)$ has both
components even. Then condition (\ref{adm3}) can be restated as no
obstacle between $O_{\zeta(\veck_{n-1})}$ and $O_{\vecl_n+\zeta(\veck_{n-1})}$, and
condition (\ref{adm4}) as the obstacle $O_{\vecl_n+\zeta(\veck_{n-1})}$ not
between $O_{\zeta(\veck_{n-1})}$ and $O_{\vecl_n+\vecl_{n+1}+\zeta(\veck_{n-1})}$.
Therefore we can take as the vertices of $G$ the pairs $(i,j)$, where
$\veci\in Q$, $\vecj\in\Z^2$, $\veci\ne \vecj$, and there is no obstacle between $O_\veci$
and $O_\vecj$. There is an edge in $G$ from $(\veci,\vecj)$ to $(\veci',\vecj')$ if and
only if $O_\vecj$ is not between $O_\veci$ and $O_{\vecj+\vecj'-\veci'}$ and $\zeta(\vecj)=\veci'$. Then,
similarly as in the torus case, there is a one-to-one correspondence
between admissible sequences and one-sided infinite paths in $G$,
starting at vertices $(\zero,\vecj)$.

This restriction on the starting point of a path in $G$ creates some
complication, but for every $\veci\in Q$ there is a similar correspondence
between admissible sequences and one-sided infinite paths in $G$,
starting at vertices $(\veci,\vecj)$. Therefore, if we want to glue finite
paths, we may choose an appropriate $\veci$.

Lemma~\ref{fingr} and its proof remain unchanged. The definition of a
\emph{small obstacle} has to be modified slightly. This is due to the
fact, that while a torus is homogeneous, so all positions of an
obstacle are equivalent, this is not the case for a square. An
obstacle placed close to a side of the square will produce a pattern
of obstacles in the unfolding which is difficult to control. Therefore
we will say that the obstacle $O$ is \emph{small} if it is contained
in a closed ball of radius smaller than $\sqrt{2}/4$, centered at
$\zero$. With this definition, Lemmas~\ref{scalar} and~\ref{small} and their proofs remain
unchanged.

We arrived at a point where the situation is completely different than
for the torus case, namely, we have to define the displacement
function. Once we do it, the definitions of the rotation set,
admissible rotation set and rotation vector remain the same, except
that we will call a rotation vector (since it belongs to $\R$) a
\emph{rotation number}.

Since we have to count how many times the trajectory rotates around
the obstacle, the simplest way is to choose a point $\vecz$ in the
interior of $O$ and set $\phi(\vecx)=\arg(\vecx-\vecz)/(2\pi)$, where arg is the
complex argument (here we identify $\R^2$ with $\C$). It is not
important that the argument is multivalued, since we are interested
only in its increment along curves. For any closed curve $\Gamma$
avoiding the interior of $O$, the increment of $\phi$ is equal to the
winding number of $\Gamma$ with respect to $\vecz$. Since the whole
interior of $O$ lies in the same component of $\C\setminus\Gamma$,
this number does not depend on the choice of $\vecz$. If $\Gamma$ is not
closed, we can extend it to a closed one, while changing the increment
of $\phi$ by less than 1. Therefore changing $\vecz$ will amount to the
change of the increment of $\phi$ by less than 2. When computing the
rotation numbers, we divide the increment of $\phi$ by the length of
the trajectory piece, and this length goes to infinity. Therefore in
the limit a different choice of $z$ will give the same result. This
proves that the rotation set we get is independent of the choice of
$\vecz$.

The proofs of the results of Section~\ref{sec_rotset} rely on the
second part of Lemma~\ref{const}, which we have not discussed yet. The
possibility of the trajectory pieces crossing the obstacles, mentioned
in that lemma, was necessary only for the proof of its first part, so
we do not have to worry about it now. However, we have to make an
additional assumption that $B$ is admissible. This creates no problem
either, because this is how we use it later. Since we changed here
many details, it makes sense to state exactly what we will be proving.

\begin{lemma}\label{new_const}
For every finite admissible sequence $B=(\veck_n)_{n=0}^s$ of elements of
$\Z^2$ the displacements $\phi$ along trajectory pieces of type $B$
differ by less than $2$.
\end{lemma}

\begin{proof}
Note that the ``folding'' map $\pi:\overline{\R^2\setminus\bigcup_{\veck
\in\Z^2}O_\veck}\to\overline{S\setminus O}$ is continuous. Therefore if
curves $\Gamma$ and $\gamma$ in $\overline{\R^2\setminus\bigcup_{\veck
\in\Z^2}O_\veck}$ with the common beginning and common end are homotopic
then $\pi(\Gamma)$ and $\pi(\gamma)$ are homotopic, so the increments
of $\phi$ along them are the same. If $\Gamma$ and $\Gamma'$ are two
trajectory pieces as in the statement of the lemma, then $\Gamma'$ can
be extended to $\gamma$ with the same beginning and end as $\Gamma$,
with the change of the increment of $\phi$ along its projection by
$\pi$ less than 2 (1 at the beginning, 1 at the end). Thus, it
suffices to show that this can be done in such a way that $\Gamma$ and
$\gamma$ are homotopic.

Therefore we have to analyze what may be the reasons for $\Gamma$ and
$\gamma$ not to be homotopic. Extending $\Gamma'$ to $\gamma$ can be
done in the right way, so this leaves two possibly bad things that we
have to exclude. The first is that when going from $O_{\veck_i}$ to
$O_{\veck_{i+2}}$ via $O_{\veck_{i+1}}$, we pass with $\Gamma$ on one side of
$O_{\veck_{i+1}}$ and with $\Gamma'$ on the other side of $O_{\veck_{i+1}}$.
The second one is that when going from $O_{\veck_i}$ to $O_{\veck_{i+1}}$, we
pass with $\Gamma$ on one side of some $O_\vecj$, and with $\Gamma'$ on
the other side of $O_\vecj$. However, the first possibility contradicts
condition (\ref{adm4}) from the definition of admissibility and the
second one contradicts condition (\ref{adm3}). This completes the
proof.
\end{proof}

With Lemma~\ref{new_const} replacing Lemma~\ref{const}, the rest of
results of Section~\ref{sec_rotset} (except the last two theorems and
the corollary) and their proofs remain the same
(with obvious minor modifications, for instance due to the fact that
constants in Lemmas~\ref{const} and~\ref{new_const} are different).
Let us state the main theorems we get in this way.

\begin{theorem}\label{new_convex}
The admissible rotation set of a billiard in a square with a small
obstacle is convex, and consequently, it is a closed interval
symmetric with respect to $0$.
\end{theorem}

\begin{theorem}\label{new_perdense}
For a billiard in a square with a small obstacle, rotation numbers of
periodic orbits of admissible type are dense in the admissible
rotation set.
\end{theorem}

\begin{theorem}\label{new_cominv}
For a billiard in a square with a small obstacle, if $u$ is a number
from the interior of $AR$, then there exists a compact, forward
invariant subset $Y$ of the phase space, such that every trajectory
from $Y$ has admissible type and rotation number $u$.
\end{theorem}

\begin{corollary}\label{new_pointwise}
For a billiard in a square with a small obstacle, if $u$ is a number
from the interior of $AR$, then there exists a trajectory of
admissible type with rotation number $u$.
\end{corollary}

\begin{corollary}\label{new_ergodic}
For a billiard in a square with a small obstacle, if $u$ is a number
from the interior of $AR$, then there exists an ergodic invariant
probability measure in the phase space, for which the integral of the
displacement is equal to $u$ and almost every trajectory is of
admissible type.
\end{corollary}

Since the billiard in the square is defined only in dimension 2, most
of the results of Section~\ref{sec_arlarge} do not have counterparts
here. However, we can investigate what happens to $AR$ as the size of
the obstacle decreases. Moreover, here the position of the obstacle
matters, so the size of the obstacle should be measured by the radius
of the smallest ball centered in the origin that contains it. The
following theorem should be considered in the context of
Theorem~\ref{wholer} that states that the full rotation set, $R$, is
equal to $[-\sqrt{2}/4,\sqrt{2}/4]$.

\begin{theorem}\label{new_epsilon}
For every $\eps>0$ there exists $\delta>0$ such that the set $AR$
contains the interval $[-\sqrt{2}/4+\eps,\sqrt{2}/4-\eps]$
whenever the obstacle is contained in the disk centered
at the origin and diameter less than $\delta$.
\end{theorem}

\begin{proof}
Let us estimate the rotation number of the curve $V$ that
in the unfolding is a straight line segment from the origin to the
point $(2n+1,2n)$. We are counting the displacement as the rotation
around the center of the square (as always, as the multiples of
$2\pi$). In particular, the displacement for $V$ makes sense, since at
its initial and terminal pieces the argument is constant.
Compare $V$ to the curve $V'$ that in the unfolding is a
straight line segment from $(0,-1/2)$ to $(2n+1,2n+1/2)$. In the
square, $V'$ goes from the lower side to the right one, to the upper
one, to the left one, to the lower one, etc., and it reflects from
each side at its midpoint. Moreover, the distances of the endpoints of
$V$ from the corresponding endpoints of $V'$ are $1/2$. Therefore,
when we deform linearly $V$ to $V'$, we do not cross any point of
$\Z^2$. This means that the difference of the displacements along
those trajectories in the square is less than 2 (this difference may
occur because they end at different points). The displacement along
$V'$ is $n+1/2$, so the displacement along $V$ is between $n-3/2$ and
$n+5/2$.

The length of $V$ is between the length of $V'$ minus 1 and the length
of $V'$, that is, between $(2n+1)\sqrt{2}-1$ and $(2n+1)\sqrt{2}$.
Therefore, as $n$ goes to infinity, the rotation number of $V$ goes to
$n/(2n\sqrt{2})=\sqrt{2}/4$.

If we fix $\eps>0$ then there is $n$ such that the rotation number of
$V$ is larger than $\sqrt{2}/4-\eps/4$. Then we can choose $\delta>0$
such that if the obstacle is contained in the disk centered at the
origin and diameter less than $\delta$ then $(2n+1,2n)$ is a vertex of
the graph $G$ and any trajectory piece $T_n$ that in the unfolding is a
straight line segment from a point of $O_\zero$ to a point of
$O_{(2n+1,2n)}$ has rotation number differing from the rotation number
of $V$ by less that $\eps/4$ (when we deform linearly $V$ to get this
trajectory piece, we do not cross any point of $\Z^2$). Hence, this
rotation number is larger than $\sqrt{2}/4-\eps/2$.

Now we construct a periodic orbit of admissible type with the rotation
number differing from the rotation number of $T_n$ by less than
$\eps/2$. By Lemma~\ref{small}, there is a loop $A_n$ in $G$, passing
through $\vecv_n=(2n+1,2n)$ and at most 2 other vertices, both from $U$.
As $n$ goes to infinity, then clearly the ratios of displacements and
of lengths of $T_n$ and $A_n$ go to 1. Therefore the ratio of their
rotation numbers also goes to 1, and if $n$ is large enough, the
difference between them will be smaller than $\eps/2$.
This gives us $\delta$ such that if the obstacle is
contained in the disk centered at the origin and diameter less than
$\delta$ then the set $AR$ contains the number $v>\sqrt{2}/4-\eps$.

Since $AR$ is symmetric with respect to 0, it contains also the number
$-v$, and since it is connected, it contains the interval
$[-\sqrt{2}/4+\eps,\sqrt{2}/4-\eps]$.
\end{proof}

Theorem~\ref{st15} also has its counterpart for the billiard in the
square.

\begin{theorem}\label{new_st15}
For a billiard in a square with a small obstacle, the admissible
rotation set is contained in the open interval
$(-\sqrt{2}/4,\sqrt{2}/4)$. In particular, $AR\ne R$.
\end{theorem}

Since the proof of this theorem utilizes a construction introduced in
the proof of Theorem~\ref{thm47}, we postpone it until the end of
Section~\ref{sec_conn}.

\section{Results on the full rotation set}\label{sec_conn}

In this section we will prove several results on the full rotation set
$R$ in both cases, not only about the admissible rotation set. Some of
the proofs apply to a much more general situation than billiards, and
then we will work under fairly general assumptions.

Let $X$ be a compact metric space and let $\Phi$ be a \emph{continuous
semiflow} on $X$. That is, $\Phi:[0,\infty)\times X\to X$ is a
continuous map such that $\Phi(0,x)=x$ and $\Phi(s+t,x)=
\Phi(t,\Phi(s,x))$ for every $x\in X$, $s,t\in[0,\infty)$. We will
often write $\Phi^t(x)$ instead of $\Phi(t,x)$. Let $\xi$ be an
\emph{time-Lipschitz continuous observable cocycle} for $(X,\Phi)$
with values in $\R^m$, that is, a continuous function
$\xi:[0,\infty)\times X\to\R^m$ such that
$\xi(s+t,x)=\xi(s,\Phi^t(x))+\xi(t,x)$ and $\|\xi(t,x)\|\le Lt$ for
some constant $L$ independent of $t$ and $x$.

The \emph{rotation set} $R$ of $(X,\Phi,\xi)$ is the set of all limits
$$\lim_{n\to\infty}\frac{\xi(t_n,x_n)}{t_n},\text{\rm\ \ where\ \ }
\lim_{n\to\infty}t_n=\infty.$$
By the definition, $R$ is closed, and is contained in the closed ball
in $\R^m$, centered at the origin, of radius $L$. In particular, $R$
is compact.

\begin{theorem}\label{conn}
The rotation set $R$ of a continuous semiflow $\Phi$ on a connected
space $X$ with a time-Lipschitz continuous observable cocycle $\xi$ is
connected.
\end{theorem}

\begin{proof}
Set
$$\psi(t,x)=\frac{\xi(t,x)}{t}$$
for $t>0$ and $x\in X$. Then the function $\psi$ is continuous on the
space $(0,\infty)\times X$. For $n\ge 1$, set $K_n=\psi([n,\infty)
\times X)$. With this notation, we have
$$R=\bigcap_{n=1}^\infty\overline{K_n}.$$
The set $[n,\infty)\times X$ is connected, so $K_n$ is connected, so
$\overline{K_n}$ is connected. Moreover, $K_n$ is contained in the
closed ball in $\R^m$, centered at the origin, of radius $L$.
Therefore $\overline{K_n}$ is compact. Thus, $(\overline{K_n})_{n=1}
^\infty$ is a descending sequence of compact connected sets, and so
its intersection $R$ is also connected.
\end{proof}

In the case of a billiard on a torus that we are considering, the
phase space $X$
is the product of the torus minus the interior of the obstacles with
the unit sphere in $\R^m$ (velocities). At the boundaries of the
obstacles we glue together the pre-collision and post-collision
velocity vectors. This space is compact, connected, and our
semiflow (even a flow, since we can move backwards in time, too) is
continuous. The observable cocycle is the displacement function.
Clearly, it is time-Lipschitz with the constant $L=1$ and continuous.
Thus, Theorem~\ref{conn} applies, and the rotation set $R$ is
connected.

Similar situation occurs in the square. Here there is one more
complication, due to the fact that there are trajectories passing
through vertices. The gluing rule at a vertex $q$ of the square is
that the phase points $(q,v)$ and $(q,-v)$ must be identified for all
relevant velocities $v$. Then the flow is also continuous in this
case, so Theorem~\ref{conn} also applies. It means that the rotation
set $R$ is a closed interval, symmetric with respect to 0.

When we work with invariant measures, we have to use a slightly
different formalism. Namely, the observable cocycle $\xi$ has to be
the integral of the \emph{observable function} $\zeta$ along an orbit
piece. That is, $\zeta:X\to\R^m$ is a bounded Borel function,
integrable along the orbits, and
$$\xi(t,x)=\int_0^t \zeta(\Phi^s(x))\;ds.$$
Assume that $\Phi$ is a continuous flow. Then, if $\mu$ is a
probability measure, invariant and ergodic with respect to $\Phi$,
then by the Ergodic Theorem, for $\mu$-almost every point $x\in X$ the
\emph{rotation vector}
$$\lim_{t\to\infty}\frac{\xi(t,x)}{t}$$
of $x$ exists and is equal to $\int_X\zeta(x)\;d\mu(x)$.

Problems may arise if we want to use weak-* convergence of measures.
If $\zeta$ is continuous and $\mu_n$ weak-* converge to $\mu$ then
the integrals of $\zeta$ with respect to $\mu_n$ converge to the
integral of $\zeta$ with respect to $\mu$. However, for a general
$\zeta$ this is not true. Note that in the cases of billiards that we
are considering, $\zeta$ is the velocity vector. It has a
discontinuity at every point where a reflection occurs (formally
speaking, it is even not well defined at those points; for definiteness
we may define it there in any way so that it remains bounded and
Borel). However, it is
well known that the convergence of integrals still holds if the set of
discontinuity points of $\zeta$ has $\mu$-measure zero (as a random
reference, we can give \cite{bauer}, Theorem~7.7.10, page 234).

Let us call an observable \emph{almost continuous} if the set of its
discontinuity points has measure zero for every $\Phi$-invariant
probability measure. By what we said above, the following lemma holds.

\begin{lemma}\label{alcont}
If probability measures $\mu_n$ weak-* converge to a $\Phi$-invar\-iant
probability measure $\mu$ and $\zeta$ is almost continuous then
$$\lim_{n\to\infty}\int_X\zeta(x)\;d\mu_n(x)=\int_X\zeta(x)\;d\mu(x).$$
\end{lemma}

We have to show that this lemma is relevant for billiards.

\begin{lemma}\label{acobs}
Let $(\Phi,X)$ be a billiard flow in the phase space. Then the
velocity observable function $\zeta$ is almost continuous,
\end{lemma}

\begin{proof}
The only points of the discontinuity points of $\zeta$ are on the
boundary of the region $\Omega$ in which we consider the billiard.
Take a small piece $Y$ of this set. Then for a small $t\ge 0$ the sets
$\Phi(t,Y)$ will be pairwise disjoint (if for $y_1\in\Phi(t_1,Y)$ and
$y_2\in\Phi(t_2,Y)$ with $t_1\ne t_2$ the velocity is the same, the
points of $\Omega$ are different). However, by the invariance of
$\mu$, their measures are the same. Since the parameter $t$ varies in
an uncountable set, those measures must be 0.
\end{proof}

The following theorem is an analogue of Theorem~2.4 from \cite{MZ2}.
Its proof is basically the same as there (except that here we deal
with a flow and a discontinuous observable), so we will omit some
details. A point of a convex set $A$ is an \emph{extreme} point of $A$
if it is not an interior point of any straight line segment contained
in $A$.

\begin{theorem}\label{extreme}
Let $(\Phi,X)$ be a continuous flow and let $\zeta:X\to\R^m$ be an almost
continuous observable function. Let $R$ be the rotation set of
$(\Phi,X,\zeta)$. Then for any extreme vector $\vecu$ of the convex hull
of $R$ there is a $\Phi$-invariant ergodic probability measure $\mu$
such that $\int_X\zeta(x)\;d\mu(x)=\vecu$.
\end{theorem}

\begin{proof}
There is a sequence of trajectory pieces such that the average
displacements on those pieces converge to $\vecu$. We can find a
subsequence of this sequence such that the measures equidistributed on
those pieces weakly-* converge to some probability measure $\nu$. This
measure is automatically invariant. Therefore, by Lemma~\ref{alcont},
we can pass to the limit with the integrals of $\zeta$, and we get
$\int_X\zeta(x)\;d\nu(x)=\vecu$. We decompose $\nu$ into ergodic
components, and since $u$ is an extreme point of the convex hull of
$R$, for almost all ergodic components $\mu$ of $\nu$ we have
$\int_X\zeta(x)\;d\mu(x)=\vecu$.
\end{proof}

By Lemma~\ref{acobs}, we can apply the above theorem to our billiards.
In particular, by the Ergodic Theorem, for any extreme vector $\vecu$ of
the convex hull of $R$ there is a point with the rotation vector $\vecu$.

Now we will look closer at billiards on the torus $\T^m$ with one
obstacle (not necessarily small). We know that its rotation set $R$ is
contained in the unit ball. It turns out that although (by
Corollary~\ref{st14}) $R$ can fill up almost the whole unit ball,
still it cannot reach the unit sphere $S^{m-1}$ on a big set. Let us
start with the following theorem.

\begin{theorem}\label{straight}
For a billiard on the torus $\T^m$ with one obstacle, if $\vecu$ is a
rotation vector of norm $1$, then there is a full trajectory in $\R^m$
which is a straight line of direction $\vecu$.
\end{theorem}

\begin{proof}
Clearly, such $\vecu$ is an extreme point of the convex hull of $R$.
By Theorem~\ref{extreme}, there is an ergodic measure $\mu$ such
that the integral of the velocity with respect to $\mu$ is $\vecu$.
Thus, the support of $\mu$ is contained in the set of points of the
phase space for which the vector component is $\vecu$. Take $t$
which is smaller than the distance between any two obstacles in the
lifting. Then $\mu$-almost all full trajectories of the billiard
have direction $\vecu$ at all times $k t$, for any integer $k$. Such
a trajectory has to be a straight line with direction $\vecu$.
\end{proof}

The following lemma has been proved in \cite{Sz} as Lemma~A.2.2.

\begin{lemma}\label{szasz}
For every dimension $m>1$ and every number $r>0$ there are finitely
many nonzero vectors $\vecx_1,\vecx_2,\dots,\vecx_k\in\Z^m$ such that whenever a
straight line $L$ in $\R^m$ is at least at the distance $r$ from
$\Z^m$, then $L$ is orthogonal to at least one of the vectors $\vecx_i$.
In other words, $L$ is parallel to the orthocomplement (lattice)
hyperplane $H_i=(\vecx_i)^\perp$.
\end{lemma}

As an immediate consequence of Theorem~\ref{straight} and
Lemma~\ref{szasz}, we get the following result.

\begin{theorem}\label{sphere}
For a billiard on the torus $\T^m$ with one obstacle, the intersection
$R\cap S^{m-1}$ is contained in the union of finitely many great
hyperspheres of $S^{m-1}$. The hyperplanes defining these
great hyperspheres can be taken as in Lemma~\ref{szasz}.
\end{theorem}

For the billiard in the square with one obstacle, we can determine the
full rotation set much better than for the torus case. In the theorem
below we do not need to assume that the obstacle is small or even
convex. However, we assume that it is contained in the interior of the
square and that its boundary is smooth.

\begin{theorem}\label{thm47}
For a billiard with one obstacle in a square, the rotation set is
contained in the interval $[-\sqrt{2}/4,\sqrt{2}/4]$.
\end{theorem}

\begin{proof}
Let $Y$ be the square minus the obstacle. In order to measure the
displacement along a trajectory piece $T$, we have to trace how its
lifting to the universal covering space of $Y$ behaves. Since $Y$ is
homeomorphic to an annulus, this universal covering has a natural
structure of a strip in the plane. Without any loss of generality, we
may assume that the displacement along $T$ is positive.

We divide $Y$ into 4
regions, as in Figure~\ref{fig47}. The line dividing regions 1
and 2 is a segment of the lowest horizontal straight line such that
the whole obstacle is below it; this segment has only its left
endpoint belonging to the obstacle. The other three dividing lines are
chosen in the same way after turning the whole picture by 90, 180 and
270 degrees. Note that here we are not interested very much to which
region the points of the division lines belong, so it is a partition
modulo its boundaries.

\begin{figure}\refstepcounter{figure}\label{fig47}\addtocounter{figure}{-1}
\begin{center}
\includegraphics[width=2.5truein]{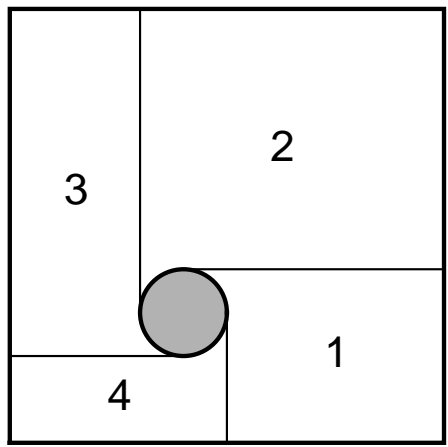}
\caption{Four regions}
\end{center}
\end{figure}

In the universal covering of $Y$, our four regions become
infinitely many ones, and they are ordered as
\begin{equation}\label{order}
\dots,1_{-1},2_{-1},3_{-1},4_{-1},1_0,2_0,3_0,4_0,1_1,2_1,3_1,4_1,
\dots
\end{equation}
Here the main number shows to which region in $Y$ our region in the
lifting projects, and the subscript indicates the branch of the
argument (as in the definition of the displacement) that we are using.
Thus, if in the lifting the trajectory goes from, say, $2_0$ to
$1_{23}$, then the displacement is 22, up to an additive constant.
This constant does not depend on the trajectory piece we consider, so
it disappears when we pass to the limit to determine the rotation set.

Since we assumed that the displacement along $T$ is positive, then in
general, the trajectory moves in the order as in (\ref{order}),
although of course it can go back and forth. Look at some region,
say, $1_n$, which is (in the universal covering) between the region
where $T$ begins and the region where $T$ ends. After $T$ leaves
$1_n$, for good, it can bounce between the left and right sides
several times. As this happens, the $y$-coordinate on the trajectory
must grow, so at some point the trajectory will hit the upper side of
the square for the first time after it leaves $1_n$ for good (unless
it ends before it does this). We denote the time of this collision by
$t(1_n)$. Then we use analogous notation for the time of the first hit
of the left side after leaving $2_n$ for good, etc.

After the trajectory hits the top side at $t(1_n)$, it moves to the
left or right (we mean the horizontal component of the velocity). It
cannot move vertically, because then it would return to the region
$1_n$. If it is moving to the left, it is still in the region $2_n$,
or it just left it, but did not hit the left side of the square yet.
Therefore $t(1_n)<t(2_n)$. If it is moving to the right, it is in, or
it will return to, the region $2_n$. Therefore also in this case
$t(1_n)<t(2_n)$.

In such a way we get an increasing sequence of times when the
trajectory $T$ hits the consecutive sides of the square (in the
lifting). By joining those consecutive reflection points by segments,
we get a piecewise linear curve $\gamma$, which is not longer than $T$,
but the displacement along $\gamma$ differs from the displacement
along $T$ at most by a constant independent of the choice of $T$. This
curve $\gamma$ goes from the right side of the square to the upper
one, to the left one, to the lower one, to the right one, etc. This is
exactly the same behavior that is displayed by the trajectory $\Gamma$
in the square without an obstacle, that starts at the midpoint of the
right side and goes in the direction of $(-1,1)$. Therefore $\gamma$
and $\Gamma$ pass through the same squares in the unfolding. We
terminate $\Gamma$ in that square in the unfolding in which $\gamma$
ends. Since in the unfolding $\Gamma$ is a segment of a straight line,
it is shorter than $\gamma$, again up to a constant independent of
$\gamma$, and those two curves have the same displacement (as always,
up to a constant). This shows that the rotation number of the
trajectory piece $T$ is not larger than for the curve $\Gamma$, plus a
constant that goes to 0 as the length of $T$ goes to infinity. Since
the rotation number of $\Gamma$ (let us think now about the infinite
trajectory) is $\sqrt{2}/4$, the limit in the definition of the
rotation set cannot exceed this number.
\end{proof}

The two trajectories in the square without an obstacle, described in
the last paragraph of the above proof, are also trajectories of any
square billiard with one small obstacle. Therefore in this case the
rotation set $R$ contains the interval $[-\sqrt{2}/4,\sqrt{2}/4]$. In
such a way we get the following result.

\begin{theorem}\label{wholer}
For a billiard with one small obstacle in a square, the rotation set
is equal to the interval $[-\sqrt{2}/4,\sqrt{2}/4]$.
\end{theorem}

Now we present the proof of Theorem~\ref{new_st15} that has been
postponed until now.

\begin{proof}[Proof of Theorem~\ref{new_st15}]
Let us use the construction from the proof of Theorem~\ref{thm47} and look
at a long trajectory piece $T$ of admissible type. We get an
increasing sequence of times when the trajectory $T$ hits the
consecutive sides of the square (in the lifting). Construct a
partial unfolding of $T$, passing to a neighboring square only at
those times. In such a way we get a piecewise linear curve $T'$ which
sometimes reflects from the sides of the unfolded square, and
sometimes goes through them. Its length is the same as the length of
$T$. Moreover, the curve $\gamma$, constructed in the proof of
Theorem~\ref{thm47} starts and terminates in the same squares as $T'$
and as we know from that proof, the displacements along $\gamma$ and
$T$ differ at most by a constant independent of $T$. It is also clear
that the same holds if we replace $\gamma$ by the segment $\Gamma$
(also from the proof of Theorem~\ref{thm47}). Thus, up to a constant
that goes to 0 as the length of $T$ goes to infinity, the rotation
number of $T$ is $\sqrt{2}/4$ multiplied by the length of $\Gamma$ and
divided by the length of $T'$.

By the same reasons as in the proof of Theorem~\ref{st15}, there exist
positive constants $c_1<c_2$ and $\alpha>0$, such that for every
trajectory piece of admissible type the distance between two
consecutive reflections from obstacles is contained in $[c_1,c_2]$, and the direction
of a trajectory piece of admissible type changes by at least $\alpha$
at each reflection. However, we have to take into account that the
direction of $T'$ can change also at the reflections from the
boundaries of the squares, and then we do not know how the angle changes.
Thus, either immediately before or immediately after each reflection from
an obstacle there must be a piece of $T'$ where the direction differs from
the direction of $\Gamma$ by at least $\alpha/2$. Such a piece has length
larger than or equal to the distance from the obstacle to the boundary
of the square, which is at least $c_1'=1/2-\sqrt{2}/4$. Those pieces are
alternating with the pieces of $T'$ of length at most $c_2$ each, where we
do not know what happens with the direction.
This means that as we follow $T'$ (except the
initial piece of the length bounded by $c_2$),
we move in the direction that differs from the direction of $\Gamma$
by at least $\alpha/2$ for at least $c_1'/(c_1'+c_2)$ of time.
Therefore, in the limit as the length of $T$ goes to infinity, the
length of $\Gamma$ divided by the length of $T'$ is not larger than
$$b=\frac{c_2}{c_1'+c_2}+\frac{c_1'}{c_1'+c_2}\cos\frac{\alpha}{2}.$$
This proves that $AR\subset [-b\sqrt{2}/4,b\sqrt{2}/4]$. Since $b<1$,
this completes the proof.
\end{proof}

\end{document}